\documentclass{article}[11pt]
\usepackage{amsthm,amsmath}
\usepackage{amssymb}
\input xypic
\usepackage[all,tips]{xy}
\newtheorem{theo}[subsubsection]{Theorem}
\newtheorem{lem}[subsubsection]{Lemma}
\newtheorem{prop}[subsubsection]{Proposition}
\newtheorem{corr}[subsubsection]{Corollary}
\newtheorem{exam}[subsubsection]{Example}
\newtheorem{rem}[subsubsection]{Remark}

\newcommand{\cC}{{\cal C}}
\newcommand{\cZ}{{\cal Z}}
\newcommand{\caD}{{\cal D}}
\newcommand{\caS}{{\cal S}}
\newcommand{\Z}{{\bf Z}}
\newcommand{\C}{{\bf C}}

\newcommand{\Rep}{{\cal R}{\it ep}}
\newcommand{\Vect}{{\cal V}{\it ect}}
\newcommand{\Hom}{{\bf Hom}}
\newcommand{\sve}{{\scriptscriptstyle {\vee}}}

\begin{document}
\author{Alexei Davydov}
\title{Modular invariants for group-theoretical modular data. I}
\maketitle
\date{}
\begin{center}
Department of Mathematics, Division of Information and Communication Sciences, Macquarie University, Sydney, NSW 2109,
Australia
\end{center}
\begin{center}
davydov@math.mq.edu.au
\end{center}
\begin{abstract}
We classify indecomposable commutative separable (special Frobenius) algebras and their local modules in (untwisted)
group-theoretical modular categories. This gives a description of modular invariants for group-theoretical modular
data. As a bi-product we provide an answer to the question when (and in how many ways) two group-theoretical modular
categories are equivalent as ribbon categories.
\end{abstract}
\tableofcontents
\section{Introduction}
An important feature of a Rational Conformal Field Theory (RCFT) is a decomposition of its partition function
$$Z(q) = \sum_{i,j}m_{i,j}\chi_i(q)\chi_j(\overline q),$$ which reflects a decomposition of the state space into a finite
sum of irreducible modules over the left-right chiral algebras. Modular invariance of the partition function implies
that the matrix of non-negative integers $M =(m_{i,j})$ is invariant with respect to the modular group actions on the
characters ({\em modular invariant}). Modules over rational chiral algebras (rational vertex operator algebras) form
modular categories \cite{ms,hu}. As an object of the category of representations of the product of the left-right
chiral algebras, the state space has a structure of commutative separable algebra \cite{sfr,hk,hk1}. Thus the problem
of classifying modular invariants (or full RCFTs) reduces to the classification of certain commutative separable
algebras in a modular category (see also \cite{os}).

One of the simplest examples of modular categories are (the categories of representations of) so-called quantum doubles
of finite groups \cite{dp}, also known as (untwisted) group-theoretical modular categories. Appearing in conformal
field theory as the modular data of holomorphic orbifolds \cite{dv,ki}, the group-theoretical modular data and
corresponding modular invariants were studied extensively (see for example \cite{cg,ev}). Relatively recently V. Ostrik
classified module categories over group-theoretical modular categories \cite{os1}, which theoretically should give the
classification of modular invariants in the case when left and right chiral modular categories coincide. The method,
used in \cite{os1}, is based on the theory of Morita equivalences for monoidal categories, developed by M. M\" uger.
Being very elegant it is also quite indirect, which unfortunately made it very difficult to calculate corresponding
modular invariants explicitly.

In this paper we describe modular invariants by classifying commutative separable (special Frobenius) algebras and
their local modules in group-theoretical modular categories. Algebras with trivial categories of so-called local
modules ({\em trivialising algebras}) correspond to modular invariants. In particular, we prove that trivialising
commutative algebras in the group-theoretical modular category, defined by a group $G$ correspond to pairs
$(H,\gamma)$, consisting of a subgroup $H\subset G$ and a 2-cocycle $\gamma\in Z^2(H,k^*)$ (which is in complete
agreement with the results from \cite{os1}). We then use the character theory for group-theoretical modular categories
to calculate corresponding modular invariants. It turns out that the character of the trivialising algebra,
corresponding to a pair $(H,\gamma)$, has the following simple form:
$$\chi(f,g) = \frac{1}{|H|}\sum_{x\in G,\\ xfx^{-1},xgx^{-1}\in H}\frac{\gamma(xfx^{-1},xgx^{-1})}{\gamma(xgx^{-1},xfx^{-1})},$$
where $f,g$ are commuting elements of $G$. By decomposing the character into a sum of irreducible characters one can
get the corresponding modular invariant. We also study trivialising algebras in a product of two group-theoretical
modular categories, corresponding to permutation modular invariants. As a result we were able to answer the question
when (and in how many ways) two group-theoretical modular categories are equivalent as ribbon categories (see also
\cite{nn}).

The paper is organized as follows. We start by listing some basic facts from the theory of modular categories, general
theory of algebras in modular categories and their relations to modular invariants (section \ref{algbmc}). Then we
study commutative separable algebras in group-theoretical modular categories (section \ref{alggrth}). We finish with
the description of modular invariants for group-theoretical modular data (section \ref{migrth}). The case of
group-theoretical modular data for the symmetric group $S_3$ is treated as an example.

\section*{Acknowledgment}
The paper was finished during the author's visit to the Max-Planck Institut f\"ur Mathematik (Bonn). The author would
like to thank MPI for hospitality and excellent working conditions. The work on the paper was partially supported by
Australian Research Council grant DP00663514. The author would like to thank I. Runkel for stimulating discussions.
Special thanks are to R. Street for invaluable support during the work on the paper.

\section{Commutative algebras in modular categories and modular invariants.}\label{algbmc}
Here we summarise some properties of and constructions associated with separable commutative algebras in braided
monoidal categories. Then we recall the notions of modular data and modular invariants and their relations to modular
categories and commutative algebras.

Throughout the paper $k$ denotes the field of complex numbers (or any other algebraically closed field of
characteristic zero). Most of our categories will be $k$-linear (all Homs are finite dimensional $k$-vector spaces,
compositions are $K$-bilinear), semi-simple (any objects is a sum of simple objects), with finitely many simple
objects. In particular, the endomorphism algebra of a simple objects is just $k$. We will denote by $Irr(\cC)$ the set
(of representatives) of isomorphism classes of simple objects in the category $\cC$. Functors are also assumed to be
$k$-linear (effects on morphisms being $k$-linear maps). A {\em fusion category} is a semi-simple $k$-linear monoidal
category, with the $k$-linear tensor product (i.e. tensor product on morphisms is $k$-linear). We also assume that the
monoidal unit of a fusion category is simple. Since it accommodates well all examples considered in this paper, we
assume that our monoidal categories are strict (associative on the nose).

\subsection{Modular categories}

Slightly changing the definition from \cite{tu} we call a fusion category {\em modular} if it is rigid, braided, ribbon
and satisfies the non-degeneracy ({\em modularity}) condition: for isomorphism classes of simple objects, the traces of
double braiding form a non-degenerate matrix
$$\tilde S = (\tilde S_{X,Y})_{X,Y\in Irr(\cC)},\quad \tilde S_{X,Y} = tr(c_{X,Y}c_{Y,X}).$$
Here $c_{X,Y}:X\otimes Y\to Y\otimes X$ is the braiding (see \cite{tu,bk} for details).

Recall that the Deligne tensor product $\cC\boxtimes\caD$ of two fusion categories is a fusion category with simple
objects $Irr(\cC\boxtimes\caD) = Irr(\cC)\times Irr(\caD)$ and the tensor product defined by
$$(X\boxtimes Y)\otimes (Z\boxtimes W) = (X\otimes Z)\boxtimes(Y\otimes W).$$
It is straightforward to see that the Deligne tensor product of two modular categories is modular.

Let $\cC$ be a ribbon category. Following \cite{tu} define $\overline{\cC}$ to be just $\cC$ as a monoidal category
with the new braiding a ribbon twist:
$$\overline{c}_{X,Y}=c_{Y,X}^{-1},\quad \overline{\theta}_X = \theta_X^{-1}.$$
Again it is very easy to see that for a modular $\cC$, $\overline\cC$ is also modular.

Examples of modular categories are provided by monoidal centre construction \cite{js}. It was proved in \cite{mu1} that
if a fusion category $\caS$ is semi-simple and spherical, then its monoidal centre $\cZ(\caS)$ is modular (see also
\cite{bv} for more general result).

\subsection{Separable algebras and their modules}

An (associative, unital) {\em algebra} in a monoidal category $\cC$ is a triple $(A,\mu,\iota)$ consisting of an
object $A\in\cC$ together with a {\em multiplication} $\mu:A\otimes A\to A$ and a {\em unit} map $\iota:1\to A$,
satisfying {\em associativiy} $$(\mu\otimes A)\mu = (M\otimes\mu)\mu,$$ and {\em unit} $$(\iota\otimes A)\mu = I =
(M\otimes\iota)\mu$$ axioms. Where it will not cause confusion we will be talking about an algebra $A$, suppressing its
multiplication and unit maps.

A right {\em module} over an algebra $A$ is a pair $(M,\nu)$, where $M$ is an object of $\cC$ and $\nu:M\otimes A\to M$
is a morphism ({\em action map}), such that $$(\nu\otimes A)\nu = (M\otimes\mu)\nu.$$ A {\em homomorphism} of right
$A$-modules $M\to N$ is a morphism $f:M\to N$ in $\cC$ such that $$(f\otimes A)\nu_N = \nu_M f.$$

Right modules over an algebra $A\in\cC$ together with module homomorphisms form a category $\cC_A$. The forgetful
functor $\cC_A\to\cC$ has a right adjoint, which sends an object $X\in\cC$ into the {\em free} $A$-module $X\otimes A$,
with $A$-module structure defined by
$$\xymatrix{X\otimes A\otimes A \ar[r]^{I\mu} & X\otimes A.}$$
Since the action map $M\otimes A\to M$ is an epimorphism of right $A$-modules any right $A$-module is a quotient of a
free module.

An algebra $(A,\mu,\iota)$ in a rigid monoidal category $\cC$ is called {\em separable} if it is equipped with a map
$\epsilon:A\to 1$ such that the following composition is a non-degenerate pairing (denoted $e:A\otimes A\to 1$)
$$\xymatrix{A\otimes A \ar[r]^(.6)\mu & A \ar[r]^{\epsilon} &  1.}$$ Non-degeneracy of $e$ means that there is a morphism
$\kappa:1\to A\otimes A$ such that the composition
$$\xymatrix{A \ar[r]^{I\kappa} & A^{\otimes 3} \ar[r]^{eI} & A}$$ is the identity. It also implies that the similar composition
$$\xymatrix{A \ar[r]^{\kappa I} & A^{\otimes 3} \ar[r]^{Ie} & A}$$ is also the identity.

Using graphical calculus for morphisms in a (rigid) monoidal category \cite{js1} one can represent morphisms between
tensor powers of a separable algebra by graphs (one dimensional CW-complexes), whose end vertices are separated into
incomig and outgoing. For example, the multiplication map $\mu$ is represented by a trivalent graph with two incoming
and one outgoing ends, the dualtity $\epsilon$ is an interval, with both incoming ends etc. It turns out (e.g it
follows from the results of \cite{rsw}) that separability implies that we can contract loops in connected graphs with
at least one end.

For a separable algebra $A$ the adjunction
$$\xymatrix{\cC \ar@/^5pt/[r] & \cC_A \ar@/^5pt/[l] }$$
splits. The splitting of the adjuction map $M\otimes A\to M$ is given by the projector $M\otimes A\to M\otimes A$:
$$\xymatrix{M\otimes A \ar[r]^{I\epsilon} & M\otimes A^{\otimes 3} \ar[r]^{I\mu I} & M\otimes A^{\otimes 2} \ar[r]^{\nu I} & M\otimes
A.}$$

For a separable algebra $A$ the effect on morphisms $\cC_A(M,N)\to\cC(M,N)$ of the forgetful functor $\cC_A\to\cC$ has
a splitting $P:\cC(M,N)\to\cC_A(M,N)$. For $f\in\cC(M,N)$ the image $P(f)$ is defined as the composition
$$\xymatrix{M \ar[r]^{I\epsilon} & M\otimes A^{\otimes 2} \ar[r]^{\nu_MI} & M\otimes A \ar[r]^{fI} & N\otimes A \ar[r]^{\nu_N} & N.}$$
Moreover, the splitting has properties $$P(fg)=fP(g)\ P(gh)=P(g)h,\quad f,h\in Mor\cC_A, g\in Mor\cC.$$ This allows to
prove {\em Maschke's lemma} for separable algebras.
\begin{lem}
Let $A$ be a separable algebra in a semi-simple rigid monoidal category $\cC$. Then the category $\cC_A$ of right
$A$-modules in $\cC$ is also semi-simple.
\end{lem}

\subsection{Local modules over commutative algebras}

A (right) module $(M,\nu)$ over a commutative algebra $A$ is {\em local} iff the diagram
$$\xymatrix{M\otimes A  \ar[r]^\nu \ar[d]_{c_{M,A}} & M\\ A\otimes M \ar[r]^{c_{A,M}} & M\otimes A \ar[u]_\nu}$$
commutes. Denote by $\cC_A^{loc}$ the full subcategory of $\cC_A$ consisting of local modules. The following result was
established in \cite{pa}.
\begin{prop}
The category $\cC_A^{loc}$ is a monoidal subcategory of $\cC$. Moreover, the braiding in $\cC$ induces the braiding in
$\cC_A^{loc}$.
\end{prop}

The following statement was proved in \cite{ffrs}.
\begin{prop}\label{algloc}
Let $(A,m,i)$ be a commutative algebra in a braided category $\cC$. Let $(B,\mu,\iota)$ be an algebra in ${\cC}_{A}$.
Define $\overline\mu$ and $\overline\iota$ as compositions $$\xymatrix{B\otimes B \ar[r] & B\otimes_AB \ar[r]^\mu & B,
& & & 1 \ar[r]^i & A \ar[r]^\iota & B.}$$ Then $(B,\overline\mu,\overline\iota)$ is an algebra in $\cC$.
\newline
The map $\iota:A\to B$ is a homomorphism of algebras in $\cC$.
\newline
The algebra $(B,\overline\mu,\overline\iota)$ in $\cC$ is separable or commutative if and only if the algebra
$(B,\mu,\iota)$ in $_{A}{\cC}$ is such.
\newline
The functor $(\cC_A^{loc})_B^{loc}\to \cC_B^{loc}$
\begin{equation}\label{locloc}
(M,m:B\otimes_A M\to M)\mapsto (M,\overline m:B\otimes M\to B\otimes_A M\stackrel{m}{\to} M)
\end{equation}
is a braided monoidal equivalence.
\end{prop}

\begin{rem}
\end{rem}
The natural map $X\otimes Y\to X\otimes_A Y$ can be seen as a lax monoidal structure on the forgetful functor
$\cC_A^{loc}\to\cC$. The commutative diagram
$$\xymatrix{ X\otimes Y \ar[r]^{c_{X,Y}} \ar[d] & Y\otimes X \ar[d]\\
X\otimes_A Y \ar[r]^{c_{X,Y}} & Y\otimes_A X}$$ implies that this lax monoidal structure is braided.
\newline
It is known that lax monoidal functors preserve structures of algebras and modules. Braided lax monoidal functors
preserve commutative algebras and local modules. This proves a half of the proposition \ref{algloc}.

We call a separable indecomposable commutative algebra $A$ in a modular category $\cC$ {\em trivializing} if
$\cC_A^{loc}=\Vect$.

\subsection{Full centre}\label{totcent}

Details of the constructions and proofs of the results of this section can be found in \cite{ffrs,kr}.

Let $A$ be an algebra in a braided category $\cC$. Its {\em left centre} $C_l(A)$ ($C_r(A)$) is an object in $\cC$ with
a morphism into $A$, universal with respect to the following property: for any $C\to A$, such that the diagram
$$\xymatrix{ C\otimes A \ar[rr] \ar[dd]_{c_{C,A}} && A\otimes A \ar[rd]^\mu \\ && & A\\ A\otimes C \ar[rr] && A\otimes
A \ar[ru]_\mu}$$ commutes, the morphism $C\to A$ factors through a morphism $C\to C_l(A)$. Right centre $C_r(A)$ is
defined similarly. The universal property implies in particular that $C_l(A)$,$C_r(A)$ are commutative algebras in
$\cC$. Note that if $A$ is a separable indecomposable algebra then $C_l(A)$,$C_r(A)$ are images of certain idempotents
on $A$ (i.e. are direct summands of $A$).

For the next construction we need to recall the fact that for modular $\cC$ the category $\cC\boxtimes\overline\cC$
contains a distinguished separable indecomposable commutative algebra $T$ (as an object $\oplus_X X\boxtimes X$ with
the sum over isomorphism classes of simple objects in $\cC$). Now the {\em full centre} of an algebra $A\in\cC$ is
$Z(A)=C_l(A\boxtimes 1)\otimes T)$ (which also equals $C_r(1\boxtimes A)\otimes T)$).

\begin{theo}
For a separable indecomposable algebra $A$ in a modular category $\cC$ the full centre $Z(A)$ is a trivializing algebra
in $\cC\boxtimes\overline\cC$.
\newline
Moreover, the full centre construction establishes an isomorphism between the set of Morita equivalence classes of
separable indecomposable algebras in $\cC$ and isomorphism classes of trivializing algebras in
$c\C\boxtimes\overline\cC$.
\end{theo}
Here two algebras in $\cC$ are Morita equivalent if their categories of modules are equivalent as module categories
over $\cC$. Hence the theorem says that the full centre is an invariant of categories of internal modules in $\cC$
(i.e. module categories over $\cC$).

\subsection{Commutative algebras in products of braided categories and their parents}\label{par}

Let $\cC\boxtimes\caD$ be the Deligne product of two braided categories. For $X\in\cC$ the functor $X\boxtimes\
:\caD\to\cC\boxtimes\caD$ has a right adjoint $\Hom_\cC(X,\ ):\cC\boxtimes\caD\to\caD$, which  can be defined as the
composition
$$\xymatrix{\cC\boxtimes\caD \ar[rrr]^{Hom_\cC(X,\ )\boxtimes I_\caD} &&& \Vect\boxtimes\caD \ar[r] & \caD . }$$
By the definition $\Hom_\cC(X,Y\boxtimes Z) = Hom_\cC(X,Y)Z$, which allows to define a map $X\boxtimes\Hom_\cC(X,A)\to
A$. The object $\Hom_\cC(X,A)$ has a universal property: the pair $(\Hom_\cC(X,\ ),X\boxtimes\Hom_\cC(X,A)\to A)$ is
terminal among the pairs $(Y,X\boxtimes Y\to A)$, i.e. for any morphism $X\boxtimes Y\to A$ in $\cC\boxtimes\caD$ there
is a unique morphism $Y\to\Hom_\cC(X,A)$, which makes the triangle
$$\xymatrix{X\boxtimes\Hom_\cC(X,A) \ar[r] & A\\ X\boxtimes Y \ar[u] \ar[ru] }$$ commute. This, in particular, can be
used to define a functorial map
$$\Hom_\cC(X,A)\otimes\Hom_\cC(X,A)\to \Hom_\cC(X\otimes Y,A\otimes B).$$
Similarly $\Hom_\caD(Y,\ ):\cC\boxtimes\caD\to\cC$ for $Y\in\caD$.
\newline
In particular, we have braided lax monoidal functors (corresponding to monoidal units in $\cC$ and $\caD$):
$$\xymatrix{\cC &&\cC\boxtimes\caD \ar[rr]^{\Hom_\cC(1,\ )}  \ar[ll]_{\Hom_\caD(1,\ )} && \caD}$$
Now, for a commutative algebra $C$ in $\cC\boxtimes\caD$, the objects $C_l=\Hom_\caD(1,C)\in\cC$,
$C_r=\Hom_\cC(1,C)\in\caD$ have the structures of commutative algebras. We call them the {\em parents} of $C$. Note
that if $C$ is indecomposable or separable, then so are its parents $C_{l,r}$.

The following theorem is a slight generalisation of results from \cite{ffrs}.
\begin{theo}
Let $C$ be a trivialising algebra in a modular category $\cC\boxtimes\caD$. Then the functor $\Hom_{\cC_{C_l}^{loc}}(\
,C)$ induces a braided monoidal equivalence $$(\cC_{C_l}^{loc})^{op}\to \caD_{C_r}^{loc}$$ of the categories of local
modules. Moreover,
\begin{equation}\label{genc}
C = \oplus_{M\in Irr(\cC_{C_l}^{loc})} M\boxtimes\Hom_\cC(M,C),
\end{equation}
where the sum is taken over simple local $C_l$-modules in $\cC$.
\newline
Conversely, for any indecomposable separable commutative algebras $A\in\cC$, $B\in\caD$ and an equivalence of braided
monoidal categories $(\cC_A^{loc})^{op}\to \caD_B^{loc}$ there exists a maximal indecomposable separable commutative
algebra in $C\in \cC\boxtimes\caD$ such that $C_l=A, C_r=B$ and the equivalences $(\cC_A^{loc})^{op}\to \caD_B^{loc}$,
$(\cC_{C_l}^{loc})^{op}\to \caD_{C_r}^{loc}$ coincide.
\end{theo}

In particular, the parents of a full centre $Z(A)$ are (left, right) centers $C_l(A),C_r(A)$.

\subsection{Modular data and modular invariants}

Recall that the {\em modular group} is the group $SL_2(\Z)$ of determinant 1 integer $2\times 2$-matrices. It is
generated by the matrices
$$s = \left(\begin{array}{rr} 0 & -1\\ 1 & 0 \end{array}\right),\quad t = \left(\begin{array}{rr} 1 & 1\\ 0 & 1 \end{array}\right),$$
with the generating system of relations $s^4 = 1, (ts)^3 = s^2$.

Let $\cC$ be a modular category. Define
$$S = (\sqrt{dim(\cC)})^{-1}\tilde S,\quad T = diag(\theta_X).$$ The proof of the following result can be found in \cite{tu}.
\begin{theo}
Let $\cC$ be a modular category. Then the operators $S$ and $T$ define a (projective) action of the modular group
$SL_2(\Z)$ on the  complexified Grothendieck group $K_0(\cC)\otimes\C$.
\end{theo}
\begin{rem}
\end{rem}
Projectivity of the above action manifests itself by a scalar multiple appearing in the second defining relation:
$$S^4 = 1,\quad (TS)^3 = \lambda S^2.$$ Over the complex numbers it is always possible to turn it into a genuine representation, by
rescaling $T$. For the reasons of why one should not do it see \cite{tu,bk}.

An alternative approach to the modular group action was developed in \cite{ly} (see also \cite{bv}). Recall  that the
{\em coend} of a monoidal category $\cC$ is an object $C\in\cC$ with a natural collection of (action) maps $X\otimes
C\to X$, universal in the following sense: for any other object $D\in\cC$ together with a natural collection of maps
$X\otimes D\to X$ there is a morphism $D\to C$ making the diagram
$$\xymatrix{ X\otimes D \ar[rr] \ar[dr] && X\otimes C \ar[dl] \\ & X}$$
commutative. Alternatively, (in the autonomous case) the coend can be defined as a colimit $\int^X X^\sve\otimes X$,
which in the case of semi-simple $\cC$ coincides with the direct sum $\oplus_X X^\sve\otimes X$ over the isomorphism
classes of simple objects. The coend has a number of nice properties and structures, e.g. for a braided $\cC$ the coend
becomes an internal Hopf algebra. In the case of a modular category $\cC$ the coend gets equipped with a projective
action of a mapping class group of a torus with removed disk. The Hom-space $\cC(1,C)$ carries a projective action of
the mapping class group of a closed torus, i.e. the modular group action. Note finally, that for a semi-simple modular
$\cC$, $\cC(1,C)$ coincides with $K_0(\cC)\otimes\C$ as a module over the modular group. The map $K_0(\cC)\otimes\C\to
\cC(1,C)$ is the composition of the {\em character map} $K_0(\cC)\otimes\C\to End(id_\cC)^\sve$ with the natural
identification $End(id_\cC)^\sve\simeq \cC(1,C)$. Here $End(id_\cC)$ is the space of endomorphisms of the identity
functor on $\cC$ and the character map sends a class of a simple object $X$ into the function $a\mapsto\lambda$ where
$a\in End(id_\cC)$ and $a_X=\lambda I_X$.

Now we explain the relation between modular invariants and trivialising algebras. The next theorem is theorem 4.5 from
\cite{ko}.
\begin{theo}
Let $A$ be an indecomposable separable commutative algebra in a modular category $\cC$ with $\theta_A = 1$. Then
$\cC_A^{loc}$ is a modular category and the map $K_0(\cC_A^{loc})\otimes_\Z\C\to K_0(\cC)\otimes_\Z\C$, induced by the
forgetful functor $\cC_A^{loc}\to\cC$ is $SL_2(\Z)$-equivariant.
\end{theo}

\begin{corr}
Let $Z$ be a trivialising algebra in a modular category $\cC$. Then its class $[A]$ in the Grothendieck ring
$K_0(\cC)\otimes_\Z\C$ is a modular invariant element.
\end{corr}
\begin{proof}
Since $A$ is a trivialising algebra, the Grothendieck group $K_0(\cC_A^{loc})$ is isomorphic to $\Z$ and the
homomorphism $K_0(\cC_A^{loc})\to K_0(\cC)$ sends an integer $n$ into $n[A]$. Modular invariance of the
complexification of this homomorphism implies that $[A]$ is a modular invariant element.
\end{proof}

It was shown in \cite{frs,hk}, that rational conformal field theories correspond to trivialising algebras in
$\cC_l\boxtimes\overline\cC_r$. Here $\cC_{l,r}$ are {\em chiral} modular categories of the theory (representation
categories of chiral vertex operator algebras). In particular, the coefficients of the decomposition of the partition
function of the theory into the sum of chiral irreducible characters are the decomposition coefficients of the
trivialising algebra in the basis of simple objects in $K_0(\cC_l\boxtimes\overline\cC_r)=K_0(\cC_l)\otimes
K_0(\cC_r)$. Traditionally \cite{ms} elements in $K_0(\cC_l)\otimes K_0(\cC_r)$, invariant with respect to the
(anti-)diagonal modular group action, are called {\em modular invariants}. A modular invariant is {\em physical} if it
corresponds to a rational conformal field theory, i.e. is the class of a trivialising algebra. In the case when
$\cC_l=\cC_r$ ({\em non-heterotic case}) are the {\em diagonal modular invariant} $\oplus_X[X]\otimes [X]$ and the {\em
conjugation modular invariant} $\oplus_X[X]\otimes [X^\sve]$. Here sums are over isomorphism classes of simple objects
of $\cC=\cC_l=\cC_r$. While the diagonal modular invariant is always physical (is the class of the full centre
$Z(1_\cC)\in\cZ(\cC)\simeq\cC\boxtimes\overline\cC$) the conjugation modular invariant can be non-physical.

\section{Commutative algebras in group-theoretical modular categories}\label{alggrth}

\subsection{Group-theoretical modular categories}
Here we describe the monoidal centre $\cZ(\cC(G))$ of the fusion category $\cC(G)$ of $G$-graded finite dimensional
vector spaces. The results of this section are mostly well-known. We will try to give references wherever it is
possible.

An {\em compatible} $G$-action on a $G$-graded vector space $V = \oplus_{g\in G}V_g$ is a collection of automorphisms
$f:V\to V$ for each $f\in G$ such that $f(V_g) = V_{fgf^{-1}}$ and $(fg)(v) = f(g(v))$.

\begin{prop}\label{mc}
The monoidal centre $\cZ(\cC(G))$ is isomorphic, as braided monoidal category, to the category $\cZ(G)$, whose objects
are $G$-graded vector spaces $X = \oplus_{g\in G}X_g$ together with a compatible $G$-action and morphisms are graded
and action preserving homomorphisms of vector spaces. The tensor product in $\cZ(G)$ is the tensor product of
$G$-graded vector spaces with the $G$-action defined by
\begin{equation}\label{tp}
f(x\otimes y) = f(x)\otimes f(y),\quad x\in X, y\in Y.
\end{equation}
The monoidal unit is $1=1_e=k$ with trivial $G$-action.
\newline
The braiding is given by
\begin{equation}\label{br}
c_{X,Y}(x\otimes y) = f(y)\otimes x,\quad x\in X_f, y\in Y.
\end{equation}
The category $\cZ(G)$ is rigid, with dual objects $X^\sve = \oplus_f (X^\sve)_f$ given by
$$(X^\sve)_f = (X_{f^{-1}})^\sve = Hom(X_{f^{-1}},k),$$ with the action
$$g(l)(x) = l(g^{-1}(x)),\quad l\in Hom(X_{f^{-1}},k), x\in X_{gf^{-1}g^{-1}}.$$
The category $\cZ(G)$ is unitarisable with the ribbon twist
$$\theta_X(x) = f^{-1}(x),\quad x\in X_f.$$
The (unitary) trace of an endomorphism $a:X\to X$ can be written in terms of ordinary traces on vector spaces $X_g$:
$$tr(a) = \sum_{g\in G}tr_{X_g}(a_g),$$
and the (unitary) dimension of an object $X\in \cZ(G)$ is the dimension of its underlying (graded) vector space
$$dim(X) = \sum_{g\in G}dim(X_g).$$
\end{prop}
\begin{proof}
For an object $(X,{\it x})$ of the centre $\cZ(\cC(G))$ the natural isomorphism $${\it x}_V:V\otimes X\to X\otimes
V,\quad V\in \cC(G)$$ is defined by its evaluations on one-dimensional graded vector spaces. Denote by $k(f)$ such a
one-dimensional graded vector space, sitting in degree $f$. Then the isomorphism ${\it x}_{k(f)}$ can be seen as an
automorphism $f:X\to X$. The fact, that ${\it x}_{k(f)}$ preserves grading, amounts to the condition $f(X_g) =
X_{fgf^{-1}}$: $$\xymatrix{X_g = (k(f)\otimes X)_{fg} \ar[rr]^{\chi_{k(f)}} && (X\otimes k(f))_{fg} = X_{fgf^{-1}}\
.}$$ The coherence condition for ${\it x}$ is equivalent to the action axioms. The diagram, defining the second
component $\chi|\psi$ of the tensor product $(X,{\it x})\otimes(Y,{\it y}) = (X\otimes Y,{\it x}|{\it y})$, is
equivalent to the tensor product of actions (\ref{tp}).
\newline
The description of the monoidal unit in a monoidal centre corresponds to the answer for the monoidal unit in $\cZ(G)$.
\newline
Clearly, the braiding $c_{(X,{\it x}),(Y,{\it y}} = {\it x}_Y$ in the centre $\cZ(\cC(G))$ corresponds to (\ref{br}).
\newline
The answer for the dual object in $\cZ(G)$ follows from the general construction of dual objects in monoidal centers of
spherical categories (see \cite{mu1}). In our concrete case it can also be verified directly. Indeed, the evaluation
map $ev_X:X^\sve\otimes X\to 1$ pairs $(X^\sve)_f$ with $X_{f^{-1}}$ via $ev_X(l\otimes x)=l(x)$. Its $G$-invariance
follows form the definition of the $G$-action on $X^\sve$:
$$ev_X(g(l\otimes x))=ev_X(g(l)\otimes g(x))=g(l)(g(x))=l(g^{-1}(g(x)))=l(x).$$
The coevaluation map $\kappa_X:1\to X\otimes X^\sve$ is defined as follows: projected to $X_g\otimes
(X^\sve)_{g^{-1}}=X_g\otimes X^*_g$ it coincides with coevaluation $\kappa_{X_g}$. The duality axioms are
straightforward.
\newline
Note that the inverse to $\theta$ has the form $\theta^{-1}(x)=f(x)$. Indeed,
$$\theta^{-1}\theta(x) = \theta^{-1}(f^{-1}(x)) =  f(f^{-1}(x)) = x.$$
The balancing axiom for $\theta$ can be checked directly. Indeed, the effect of the double braiding on $x\otimes y\in
X_f\otimes Y_g$ is
$$x\otimes y\mapsto f(y)\otimes x\mapsto (fgf^{-1}(x)\otimes f(y),$$ while $\theta_{X\otimes
Y}^{-1}(\theta_X\otimes\theta_Y)$ acts as
$$x\otimes y\mapsto f^{-1}(x)\otimes g^{-1}(y)\mapsto (fg)(f^{-1}(x)\otimes g^{-1}(y)) =$$
$$(fg)(f^{-1}(x))\otimes (fg)(g^{-1}(y)) = (fg(f^{-1})(x)\otimes f(y).$$
The self-duality for the ribbon twist $\theta_{X^\sve}=(\theta_X^{-1})^\sve$ is straightforward.
\newline
The formula for the trace follows from the fact that the duality structure in $\cZ(G)$ is the same as in the category
of finite dimensional ($G$-graded) vector spaces.
\end{proof}

In the next statement we describe simple objects and the $S$-matrix of the category $\cZ(G)$ (see also \cite{cg}).
\begin{prop}\label{simpl}
Simple objects of $\cZ(G)$ are parametrised by pairs $(g,U)$, where $g\in G$ and $U$ is a simple module over the
twisted group algebra $k[C_G(g)]$.
\newline
The dimension of the category $\cZ(G)$ is $|G|^2$.
\newline
The category $\cZ(G)$ is modular with the $S$- and $T$-matrices:
$$S_{(f,\psi),(g,\xi)} = \frac{1}{|G|}\sum_{u\in f^G, v\in g^G, uv=vu}\psi(xv^{-1}x^{-1})\xi(yu^{-1}y^{-1}),$$ where $u=x^{-1}fx,\ v=y^{-1}gy$, and
\begin{equation}\label{tmat}
T_{(f,\psi),(f,\psi)} = \frac{\psi(f)}{\psi(e)}.
\end{equation}
\end{prop}
\begin{proof}
Clearly the support of a simple object $V$ in $\cZ(G)$ should be an indecomposable $G$-subset in $G$ (with conjugation
action), i.e. a conjugacy class of $G$. Let $g$ be an element of the support. The axioms of the action imply that $V$
is induced from the $k[C_G(g)]$-module $V_g$. Finally, for $V$ to be simple, the $k[C_G(g)]$-module $V_g$ must be
simple as well.
\newline
For $g\in G$ the sum $\sum_U dim(U)^2$ over isomorphism classes of irreducible $k[C_G(g)]$-modules is equal to
$|C_G(g)|$. Since $dim(g,U)=dim(Ind_{C_G(g)}^G(U())=[G:C_G(g)]dim(U)$ $$dim(\cZ(G)) = \sum_{g,U}[G:C_G(g)]^2dim(U)^2 =
\sum_g [G:C_G(g)]^2|C_G(g)| =$$ $$ |G|^2\sum_g|C_G(g)|^{-1},$$ where $g$ runs through representatives of conjugacy classes
of $G$. It is well-known in group theory that the last sum is equal to $|G|^2$.
\newline
The formula for the $S$-matrix can be obtained by calculating the trace of the double braiding
$c_{(g,\xi),(f,\psi)}c_{(f,\psi),(g,\xi)}$ in the category $\cZ(G)$.
\end{proof}

The next result describes Deligne products and mirrors of group-theoretical modular categories.
\begin{prop}
$$\cZ(G_1)\boxtimes\cZ(G_2)\simeq \cZ(G_1\times G_2),\quad
\overline{\cZ(G)}\simeq \cZ(G).$$
\end{prop}
\begin{proof}
Follows from the straightforward equivalences:
$$\cZ(\cC)\boxtimes\cZ(\caD)\simeq\cZ(\cC\boxtimes\caD),\quad \overline\cZ(\cC)\simeq\cZ(\cC^{op}).$$
\end{proof}

For the case $\cC=\cZ(G)$ the space of characters has the following description (see \cite{ba1,ba,wi}). It is the space
of $k$-valued functions on
$$C^2_\alpha(G) = \{(f,g)\in G^{\times 2},\quad fg=gf\}.$$
In this realisation the $SL_2(\Z)$-action is given by
$$S(\chi)(f,g) = \chi(g,f^{-1}),\quad T(\chi)(f,g) = \chi(f,fg).$$

For an object $X$ the character map $K_0(\cZ(G))\otimes_\Z\C\to Hom(C^2_\alpha(G),k)$ sends the class $[X]$ into the
function (the {\em character}): $$\chi_X(f,g) = tr_{X_f}(g).$$ In particular, the character of the dual object (the
{\em dual character}) has a form:
$$\chi_{X^\sve}(f,g) = \chi_X(f^{-1},g^{-1}).$$
As in the ordinary character theory, the space of characters of $\cZ(G)$ comes equipped with a scalar product (see
\cite{ba1})
$$(\chi,\psi) = \frac{1}{|G|}\sum_{f,g\in G}\chi(f,g)\overline{\psi(f,g)},$$ which calculates dimensions of
corresponding Hom-spaces in $\cZ(G)$:
$$(\chi_X,\chi_Y) = dim(\cZ(G)(X,Y)).$$ In particular, for
irreducible $X,Y$, $(\chi_X,\chi_Y) = 1$ iff $X=Y$ and zero otherwise.

\subsection{Algebras in group-theoretical modular categories}

We start with expanding the structure of an algebra in the category $\cZ(G)$ in plain algebraic terms. Recall that a
$G$-graded vector space $A = \oplus_{g\in G}A_g$ is a {\em $G$-graded algebra} if the multiplication preserves grading
$A_fA_g\subset A_{fg}$.
\begin{prop}
An algebra in the category $\cZ(G)$ is a $G$-graded associative algebra together with a $G$-action such that
\begin{equation}\label{ah}
f(ab) = f(a)f(b),\quad a,b\in A.
\end{equation}
An algebra $A$ in the category $\cZ(G)$ is commutative iff
\begin{equation}\label{co}
ab = f(b)a,\quad \forall a\in A_f, b\in A.
\end{equation}
The twist $\theta_A$ is trivial iff
$$f(a) = a,\quad a\in A_f.$$
\end{prop}
\begin{proof}
Being a morphism in the category $\cZ(G)$ the multiplication of an algebra in $\cZ(G)$ preserves grading and $G$-action
(hence the property (\ref{ah})). Associativity of multiplication in $\cZ(G)$ is equivalent to ordinary associativity.
\newline
The formula (\ref{br}) for the braiding in $\cZ(G)$ implies that commutativity for an algebra $A$ in the category
$\cZ(G)$ is equivalent to the condition (\ref{co}).
\end{proof}

By {\em $G$-algebra} we mean an algebra with an action of $G$ by algebra homomorphisms. Note that the degree $e$ part
$A_e$ of an algebra $A$ in the category $\cZ(G)$ is an associative $G$-algebra and $A$ is a module over $A_e$. Moreover
the algebra $A_e$ is commutative if $A$ is a commutative algebra in the category $\cZ(G)$.

\begin{prop}\label{sep}
An algebra $A$ in the category $\cZ(G)$ is separable iff
$$\xymatrix{A_f\otimes A_{f^{-1}} \ar[r]^\mu & A_e \ar[r]^{\epsilon} & k}$$
defines a non-degenerate bilinear pairing for any $f\in G$. In particular, the algebra $A_e$ is separable if $A$ is a
separable algebra in the category $\cZ(G)$.
\end{prop}
\begin{proof}
Being a graded homomorphism the separability map $A\to 1$ is zero on $A_f$ for $f\not= e$. Hence the separability
bilinear form is zero on $A_f\otimes A_g$ unless $fg=e$. In particular, the restriction of $\epsilon$ to $A_e$ makes it
a separable algebra in the category of vector spaces.
\end{proof}

\subsection{Commutative separable algebras in trivial degree and their local modules}
We start with a well known (see for example \cite{ko}) description of indecomposable commutative separable
$G$-algebras. We give (a sketch of) the proof for completeness.
\begin{lem}\label{comgalg}
Commutative separable $G$-algebras are function algebras on $G$-sets. Indecomposable $G$-algebras correspond to
transitive $G$-sets.
\end{lem}
\begin{proof}
A separable commutative algebra over an algebraically closed field is a function algebra $k(X)$ on a finite set $X$
(with elements of $X$ corresponding to minimal idempotents of the algebra). The $G$-action on the algebra amounts to a
$G$-action on the set $X$. Obviously, the algebra of functions $k(X\cup Y)$ on the disjoint union of $G$-sets is the
direct sum of $G$-algebras $k(X)\oplus k(Y)$ and any direct sum decomposition of $G$-algebras appears in that way.
\end{proof}

Let $k(X)$ be an indecomposable $G$-algebra. By choosing a minimal idempotent $p\in X$, we can identify the $G$-set $X$
with the set $G/H$ of cosets modulo the stabilizer subgroup $H=St_G(p)$.
\begin{theo}\label{loc1}
The category ${\cZ(G)}^{loc}_{k(G/H)}$, of local left $k(G/H)$-modules in $\cZ(G)$, is equivalent, as a ribbon
category, to $\cZ(H)$.
\end{theo}
\begin{proof}
For a right $k(G/H)$-module $M$ the product $Mp$ with a chosen idempotent is a $G$-graded vector space with $H$-action.
For a local $M$ the support of $Mp$ (elements of $G$, whose graded components are non-zero) is a subset of $H$. Indeed,
for $m\in M_f$ the locality condition implies that $mp = mf(p)$ and $mp = mp^2 = mpf(p)$. Thus if $mp\not=0$ the
product $pf(p)$ is also non-zero and $f(p)=p$. Hence for a local $H$ the subspace $Mp$ is an object of
$\cZ_{\alpha|_H}(H)$, which defines a functor $${\cZ(G)}^{loc}_{k(G/H)}\to \cZ_{\alpha|_H}(H),\quad M\mapsto Mp.$$ The
functor is obviously monoidal $(M\otimes_A N)p = Mp\otimes Np$, braided and balanced.

Now let $U\in \cZ(H)$. The tensor product $k(G)\otimes_H U$ (which is spanned by $p_g\otimes u$, modulo $p_{gh}\otimes
u = p_g\otimes h(u)$) is naturally equipped with the $G$-grading $$|p_g\otimes u| = g|u|g^{-1}$$ and the $G$-action
$f(p_g\otimes u) = p_{fg}\otimes u,$ making it an object of $\cZ(G)$. The homomorphism of algebras $k(G/H)\to k(G)$
(induced by the quotient map $G\to G/H$) makes $k(G)\otimes_H U$ a right $k(G/H)$-module. Explicitly, for a coset $x\in
G/H$ $$(p_g\otimes u)p_x = \delta_{g,x}p_g\otimes u.$$ Here $\delta_{g,x}$ is the $\delta$-function, which is equal to
1, if $g$ belongs to $x$, and zero otherwise. Moreover, $k(G)\otimes_H U$ is a local left $k(G/H)$-module: the value of
the product map on $p_x\otimes(p_g\otimes u)$ coincides with the value on
$$(c\circ c)(p_x\otimes(p_g\otimes u)) =
g|u|g^{-1}(p_x)\otimes(p_g\otimes u) = p_{g|u|g^{-1}x}\otimes(p_g\otimes u).$$ Indeed, $g$ belongs to $x$ (i.e. $x =
gH$) iff $g$ belongs to $g|u|g^{-1}x = g|u|g^{-1}gH = g|u|H = gH$. Thus we have a functor
$$\cZ(H)\to {\cZ(G)}^{loc}_{k(G/H)},\quad U\mapsto k(G)\otimes_H U.$$ Finally, the maps $$U\to
(k(G)\otimes_H U)p,\quad u\mapsto p_e\otimes u,$$ $$k(G)\otimes_H M\to M,\quad p_g\otimes mp\mapsto g(mp)$$ are
isomorphisms.
\end{proof}
\begin{rem}\end{rem}
It follows from the proof of the theorem \ref{loc1} that the category ${\cZ(G)}_{k(G/H)}$ of right $k(G/H)$ modules can
be identified with the category of $G$-graded vector spaces equipped with $H$-actions.

\begin{rem}\end{rem} Theorem \ref{loc1} in combination with proposition \ref{algloc} gives an interpretation of the {\em transfer},
defined in \cite{tu}. The transfer turns an algebra from $\cZ(H)$ into an algebra from $\cZ(G)$. Indeed, by theorem
\ref{loc1} an algebra from $\cZ(H)$ is an algebra in ${\cZ(G)}^{loc}_{k(G/H)}$, which by proposition \ref{algloc} gives
an algebra in $\cZ(G)$.

\begin{corr}\label{trans}
For a simple separable algebra $A$ in $\cZ(G)$ there is a subgroup $H\subset G$ such that $A$ is the transfer of a
simple separable algebra $B$ in $\cZ(H)$ with $B_e=k$.
\end{corr}
\begin{proof}
The subalgebra $A_e$ is an indecomposable commutative $G$-algebra. By lemma \ref{comgalg} it is isomorphic to $k(X)$
for some transitive $G$-set $X$. By proposition \ref{algloc}, $A$ is a commutative algebra in ${\cZ(G)}^{loc}_{A_e}$.
Thus, by theorem \ref{loc1}, $A$ is the transfer of the indecomposable separable algebra $B=pA$ from $\cZ(H)$ (here $p$
is the minimal idempotent of $A_e$, corresponding to an element of $X$, with the stabiliser $H=St_G(p)$). Finally, $B_e
= pA_e = k$ by minimality of $p$.
\end{proof}

\subsection{Commutative separable algebras trivial in trivial degree and their local modules}
Here we describe simple commutative separable algebras $B$ in $\cZ(H)$ with $B_e=k$.
\begin{lem}
Let $B$ be a separable algebra in $\cZ(H)$ such that $B_e=k$. Then $$dim(B_h)\leq 1,\quad \forall h\in H.$$ Moreover
the support of $B$ $$F = \{f\in H|\ B_f\not=0\}$$ is a normal subgroup of $H$.
\end{lem}
\begin{proof}
By the proposition \ref{sep} an algebra $B$, such that $B_e=k$, is separable iff the multiplication defines the
non-degenerate pairing $m:B_g\otimes B_{g^{-1}}\to A_e = k$. Thus, associativity of multiplication implies that, for
any $a,c\in B_g$ and $b\in B_{g^{-1}}$ $m(a,b)c = am(b,c).$ For non-zero $a,c$, choosing $b$ such that
$m(a,b),m(b,c)\not=0$, we get that $a$ and $c$ are proportional.
\newline
Now, it follows from the non-degeneracy of $m:B_g\otimes B_{g^{-1}}\to A_e = k$, that a generator of a non zero $B_f$
is invertible. Thus, for non-zero components $B_f,B_g$ the product $B_fB_g$ is also non-zero.
\end{proof}

Let $F\vartriangleleft H$ be a normal subgroup and $\gamma\in Z^2(F,k^*)$ be a normalised cocycle, i.e.
$\gamma(e,g)=\gamma(f,e)=1$ and
$$\gamma(f,g)\gamma(fg,h) = \gamma(g,h)\gamma(f,gh).$$
Note that for a 2-cocycle $\gamma\in Z^2(G,k^*)$ the expression
$$\gamma_f(g) = \gamma(f,g)\gamma(g,f)^{-1}$$ define a multiplicative map (character) $\gamma_f:C_G(f)\to k^*$ of the
centraliser $C_G(f)$.

Denote by $k[F,\gamma]$ an $H$-graded associative algebra with the basis $e_f, f\in F$, graded as $|e_f| = f$, and with
multiplication defined by $e_fe_g = \gamma(f,g)e_fg$.
\begin{prop}\label{ttd}
An indecomposable commutative separable algebra $B$ in $\cZ_{\alpha}(H)$ with $B_e=k$ has a form $k[F,\gamma]$ with the
$H$-action given by: $$h(e_f) = \varepsilon_h(f)e_{hfh^{-1}},$$ for some $\varepsilon:H\times F\to k^*$ satisfying
\begin{equation}\label{eps1}
\varepsilon_{gh}(f)  =  \varepsilon_g(hfh^{-1})\varepsilon_h(f),\quad g,h\in H, f\in F
\end{equation}
\begin{equation} \label{eps2}
\gamma(f,g)\varepsilon_h(fg)  =  \varepsilon_h(f)\varepsilon_h(g)\gamma(hfh^{-1}hgh^{-1}),\quad h\in H, f,g\in F
\end{equation}
\begin{equation}\label{eps3}
\gamma(f,g)  =  \varepsilon_f(g)\gamma(fgf^{-1},f),\quad f,g\in F.
\end{equation}
\end{prop}
\begin{proof}
 Indeed, action axiom requires that
$(gh)(e_f) = \varepsilon_{gh}(f)e_{ghfh^{-1}g^{-1}}$ coincides with $$g(h(e_f)) =
\varepsilon_h(f)\varepsilon_g(hfh^{-1}e_{ghfh^{-1}g^{-1}},$$ which gives the first identity. Multiplicativity of the
action amounts to the equality between $$h(e_fe_g) = \gamma(f,g)\varepsilon_h(fg)e_{hfgh^{-1}}$$ and
$$h(e_f)h(e_g) =\varepsilon_h(f)\varepsilon_h(g)\gamma(hfh^{-1},hgh^{-1})e_{hfgh^{-1}},$$
which gives the second identity. Finally, commutativity implies that $e_fe_g = \gamma(f,g)e_{fg}$ is equal to
\begin{equation}\label{coalg}
f(e_g)e_f = \varepsilon_f(g)e_{fgf^{-1}}e_f = \varepsilon_f(g)\gamma(fgf^{-1},f)e_{fg}.
\end{equation}
\end{proof}
Denote by $k[F,\gamma,\varepsilon]$ an indecomposable commutative separable algebra in $\cZ(H)$, defined in proposition
\ref{ttd}.
\begin{lem}
Two algebras $k[F,\gamma,\varepsilon]$ and $k[F',\gamma',\varepsilon']$ in the category $\cZ(H)$ are isomorphic iff
there is a cochain $c:F\to k^*$ such that
$$c(fg)\gamma(f,g) = \gamma'(f,g)c(f)c(g),\quad \varepsilon_h(f)c(hfh^{-1}) = c(f){\varepsilon '}_h(g).$$
\end{lem}
\begin{proof}
Isomorphic algebras in $\cZ(H)$ have to have the same supports. Thus $F=F'$. Since the components of both
$k[F,\gamma,\varepsilon]$ and $k[F',\gamma',\varepsilon']$ are all one dimensional, an isomorphism
$k[F,\gamma,\varepsilon]\to k[F',\gamma',\varepsilon']$ has a form $e_f\mapsto c(f)e_f$ for some $c(f)\in k^*$.
Finally, multiplicativity of this mapping is equivalent to the first condition, while $H$-equivariance is equivalent to
the second.
\end{proof}

%$C^2(F,k^*)=C^0(H,C^2(F,k^*))$. Thus, in the terminology of the appendix, $(\varepsilon,\gamma)$ is a 2-cochain of
%equivalent to the coboundary condition $d(\varepsilon,\gamma)=(\alpha_2,\alpha_1,\alpha_0)=\tau(\alpha)$ in $\tilde
%C^*(H,F,k^*)$. The equations (\ref{coalg}) say that $(\varepsilon,\gamma)=d(c)(\varepsilon',\gamma')$ for $c\in
%C^1(F,k^*)=\tilde C^1(H,F,k^*)$.

For the sake of keeping it short we will not give complete description of the category of local modules over the
algebra $k[F,\gamma,\varepsilon]$ (which will be given in the subsequent paper). Instead we characterise those algebras
which have trivial category of local modules (i.e. trivialising algebras).
\begin{theo}\label{loc2}
The algebra $k[F,\gamma,\varepsilon]$ in the category $\cZ(H)$ is trivialising iff $F=H$.
\end{theo}
\begin{proof}
The structure of a right $k[F,\gamma,\varepsilon]$-module on an object $M=\oplus_{h\in H}M_h$ of $\cZ_\alpha(H)$
amounts to a collection of isomorphisms $e_f:M_h\to M_{hf}$ (right multiplication by $e_f\in k[F,\gamma,\varepsilon]$)
such that $$e_e = I,\quad e_fe_{f'} = \gamma(f,f')e_{f'f},\quad he_fh^{-1} = \varepsilon_h(f)e_{hfh^{-1}},\quad f,f'\in
F, h\in H.$$ Here $h:M_{h'}\to M_{hh'h^{-1}}$ is the $\alpha$-projective $H$-action on $M$. The
$k[F,\gamma,\varepsilon]$-module $M$ is local iff $e_f = \varepsilon_h(f)hfh^{-1}e_{hfh^{-1}}$ on $M_h$. Indeed, the
double braiding in $\cZ(H)$ transforms an element $m\otimes e_f\in M\otimes A$ (with $m\in M_h$) as follows
$$m\otimes e_f\mapsto h(e_f)\otimes m=\varepsilon_h(f)e_{hfh^{-1}}\otimes m\mapsto\varepsilon_h(f)hfh^{-1}(m)\otimes
e_{hfh^{-1}}.$$ An equivalent way of expressing the locality condition is:
$$f=\varepsilon_h(h^{-1}fh)^{-1}\gamma(h^{-1}fh,f^{-1})\gamma(f,f^{-1})^{-1}e_{h^{-1}fhf^{-1}}=\epsilon_h(f)e_{[h^{-1},f]}.$$
In particular, $F$ acts trivially on $M_e$.
\newline
Now if $F=H$ the action map $M_e\otimes A\to M$ is an isomorphism, i.e. any local module is free. Conversely for
$F\not= H$ take a non-trivial $H/F$-representation $U$ and define an $H$-action on the $H$-graded vector space
$$M=V\otimes A = \bigoplus_{f\in F}M_f,\quad M_f=V\otimes e_f$$ by
$$h(v\otimes e_f) = \varepsilon_h(f)h(v)\otimes e_{hfh^{-1}}.$$
Then $M$ is a right $k[F,\gamma,\varepsilon]$-module
$$v\otimes e_f\otimes e_{f'}\mapsto \gamma(f,f')v\otimes e_{ff'},$$ which is local and non-free.
\end{proof}

\subsection{Commutative separable algebras and their local modules}

In this section we combine the previous results on commutative separable algebras in group-theoretical modular
categories and on their local modules.

Define $A(H,F,\gamma,\varepsilon)$ as a vector space, spanned by $a_{g,f}$, with $g\in G,f\in F$, modulo the relations
$$a_{gh,f} = \varepsilon_h(f)a_{g,hfh^{-1}},\quad \forall h\in H,$$ with the $G$-grading, given by $|a_{g,f}| =
gfg^{-1}$, the $G$-action $g'(a_{g,f}) = a_{g'g,f}$ and the multiplication
$$a_{g,f}a_{g',f'} = \delta_{g,g'}\gamma(f,f')a_{g,ff'}.$$

\begin{theo}\label{main}
Indecomposable separable commutative algebras in $\cZ(G)$ correspond to quadruples $(H,F,\gamma,\varepsilon)$, where
$H\subset G$ is a subgroup, $F\lhd H$ is a normal subgroup, $\gamma\in Z^2(F,k^*)$ is a cocycle and
$\varepsilon:H\times F\to k^*$ satisfies the conditions (\ref{eps1},\ref{eps2},\ref{eps3}).
\end{theo}
\begin{proof}
Follows from corollary \ref{trans} and propositions \ref{ttd} and \ref{algloc}.
\end{proof}

\begin{rem}
\end{rem}
Note that the twist $\theta_A$ is always trivial on the algebra $A = A(H,F,\gamma,\varepsilon)$. Indeed,
$$\theta^{-1}_A(a_{g,f}) = (gfg^{-1})(a_{g,f}) = a_{gfg^{-1}g,f} = a_{gf,f} =\varepsilon_f(f)a_{g,f}$$ with
$\varepsilon_f(f)=\gamma(f,f)\gamma(f,f)^{-1}=1$ by (\ref{eps3}).

\begin{theo}\label{local}
The algebra $A(H,F,\gamma,\varepsilon)$ in the category $\cZ(G)$ is trivialising iff $F=H$.
\end{theo}
\begin{proof}
Follows from theorems \ref{loc1},\ref{loc2} and proposition \ref{algloc}.
\end{proof}
Note that when $F=H$ the map $\varepsilon$ is completely determined by $\gamma$. Thus trivialising algebras in $\cZ(G)$
correspond to pairs $(H,\gamma)$, where $H\subset G$ is a subgroup and $\gamma\in Z^2(H,k^*)$ is a 2-cocycle.
\begin{rem}
\end{rem}
It follows from the theorem that trivialising algebras in $\cZ(G)\boxtimes \cZ(G)\simeq \cZ(G\times G)$ correspond to
pairs $(U,\gamma)$, where $U\subset G\times G$ is a subgroup and $\gamma\in Z^2(U,k^*)$ is a 2-cocycle. This coincides
with the parametrisation of module categories obtained in \cite{os}, which illustrates the fact (formulated in section
\ref{totcent}) that the total centre defines a bijection between equivalence classes of indecomposable module
categories over $\cZ(G)$ and maximal indecomposable separable commutative algebras in $\cZ(G)\boxtimes \cZ(G)$.

\subsection{Trivialising algebras in products of group-theoretical module categories and equivalences between group-theoretical
module categories}

In this section we describe the parents of maximal indecomposable commutative separable algebras in
$\cZ(G)\boxtimes\cZ(Q)$ and use this description to analyze braided monoidal equivalences between $\cZ(G)$ and
$\cZ(Q)$.

Let $G,Q$ be finite groups. It is straightforward to see that the functor, defined in section \ref{par},
$$Hom_{\cZ(G)}(1,\ ):\cZ(G\times Q) \simeq \cZ(G)\boxtimes\cZ(Q)\to\cZ(Q)$$
sends $X$ into subspace of invariants $(\oplus_{q\in Q}X_{(e,q)})^{G\times\{e\}}$. Let $U\subset G\times Q$ be a
subgroup and $\gamma\in Z^2(G\times Q,k^*)$ be a normalised 2-cocycle. The pair $(U,\gamma)$ defines a maximal
indecomposable commutative separable algebra $A(U,\gamma)$ in $\cZ(G\times Q) \simeq \cZ(G)\boxtimes\cZ(Q)$.
\begin{theo}\label{pare}
The parent $Hom_{\cZ(G)}(1,A)\in \cZ(Q)$ of the maximal indecomposable commutative separable algebra $A(U,\gamma)$ in
$\cZ(G)\boxtimes\cZ(Q)$ is isomorphic to $A(pr_2(U),K,\gamma|_K,\varepsilon)$, where $pr_2(U)\subset Q$ is the
projection of $U\subset G\times Q$ onto the second factor, $K$ is the kernel of the homomorphism
$$\overline\gamma:U\cap(\{e\}\times Q)\to \widehat{U\cap(G\times\{e\})},\quad (e,q)\mapsto \gamma_{(e,q)},$$ and
$\varepsilon:pr_2(U)\times K\to k^*$ is given by
$$\varepsilon_q(v) = \gamma((g,q)|v) = \gamma((g,q),v)\gamma((g,q)v(g,q)^{-1},(g,q)),\quad q\in pr_2(U), v\in K.$$
\end{theo}
\begin{proof}
As was noted in section \ref{par}, the algebra $B=Hom_{\cZ(G)}(1,A)$ is an indecomposable commutative separable algebra
in $\cZ_\beta(Q)$. Thus, by theorem \ref{main}, it should have a form $A(H,F,\gamma,\varepsilon)$ for some $F\lhd
H\subset Q$. To find $H$ we need to look at the trivial degree component $B_e$. Since
$$B_e = A_{(e,e)}^{G\times\{e\}} = k(G\times Q/U)^{G\times\{e\}} = k((G\times\{e\})\setminus G\times Q/U)$$ $H$ can be
defined as the stabiliser of $(e,e)$ with respect to the (transitive) $Q$-action on $(G\times\{e\})\setminus G\times
Q/U$, which coincides with
$$\{q\in Q|\ \exists g\in G: (g,q)\in U\} = pr_2(U).$$
Thus as a $Q$-algebra $B_e=k(Q/pr_2(U)$. To determine the rest of the defining data for $B$ we need to look at $pB$,
where $p$ is a minimal idempotent of $B_e$. Let $p\in B_e$ be the minimal idempotent, corresponding to the unit element
$e\in Q$. As an element of $A_{(e,e)}$ it has the following decomposition $p=\sum_{g\in G}g(\tilde p)$, where $\tilde
p$ is the minimal idempotent in $A_{(e,e)}$ corresponding to $(e,e)$. Hence
$$pB = (\sum_{g\in G}g(\tilde p))(\oplus_{q\in Q}A_{(e,q)})^{G\times\{e\}} = (\oplus_{q\in Q}\tilde pA_{(e,q)})^{(G\times\{e\})\cap U}.$$
Now, since $\tilde pA = k[U,\gamma]$, we have that $\oplus_{q\in Q}\tilde pA_{(e,q)} = k[U\cap(\{e\}\times Q),\gamma]$.
The conjugation action of $U\cap(G\times\{e\})$ on $U\cap(\{e\}\times Q)$ is trivial, so the only non-triviality comes
from $\gamma$: for $u\in U\cap(G\times\{e\}), v\in U\cap(\{e\}\times Q)$
$$u(e_v) = e_ue_ve_u^{-1} = \gamma(u|v)e_v,\quad \gamma(u|v) = \gamma(u,v)\gamma(v,u)^{-1}.$$
Note that, restricted to $(U\cap(G\times\{e\}))\times(U\cap(\{e\}\times Q))$, $\gamma(\ |\ )$ is a bi-multiplicative
pairing. Hence $e_v$ is an invariant iff $\gamma(\ |v)$ is trivial. Thus $$pB = k[U\cap(\{e\}\times
Q),\gamma]^{(G\times\{e\})\cap U} = k[K,\gamma|_K].$$ Finally, to determine $\varepsilon:pr_2(U)\times K\to k^*$ we
need to write the conjugation action of $pr_2(U)$ on $k[K,\gamma|_K]$ in the form
$(g,q)(e_v)=\varepsilon_q(v)e_{qvq^{-1}}$ for $q\in pr_2(U), v\in K$. Since
$$(g,q)(e_v) = e_{(g,q)}e_ve_{(g,q)}^{-1} = \gamma((g,q),v)\gamma((g,q)v(g,q)^{-1},(g,q))e_{qvq^{-1}}$$
we have the description for $\varepsilon$. Note that for $q\in pr_2(U), v\in K$ the value of $\gamma((g,q)|v)$ does not
depend on the choice of $g$.
\end{proof}
\begin{rem}
\end{rem}
Similarly, the parent $Hom_{\cZ(Q)}(1,A)\in\cZ(G)$ of $A(U,U,\gamma)\in \cZ(G)\boxtimes\cZ(Q)$ is isomorphic to
$A(pr_1(U),K,\gamma|_K,\varepsilon)$, where $pr_1(U)\subset G$ is the projection of $U\subset G\times Q$ onto the first
factor, $K$ is the kernel of the homomorphism $\overline\gamma:U\cap(G\times\{e\})\to \widehat{U\cap(\{e\}\times Q)}$,
induced by $\gamma(\ |\ )$, and $\varepsilon:pr_1(U)\times K\to k^*$ is given by
$$\varepsilon_g(v) = \gamma((g,q)|v) = \gamma((g,q),v)\gamma((g,q)v(g,q)^{-1},(g,q)),\quad g\in pr_1(U), v\in K.$$

\begin{corr}\label{eqbc}
An equivalence between $\cZ(G)$ and $\cZ(Q)$, as ribbon categories, corresponds to a subgroup $U\subset G\times Q$,
such that $pr_1(U) = G, pr_2(U) = Q$, together with a 2-cocycle $\gamma\in Z^2(U,k^*)$, such that $\gamma(\ |\ )$
induces a non-degenerate pairing
\begin{equation}\label{pair}
(U\cap(G\times\{e\}))\times(U\cap(\{e\}\times Q))\to k^*.
\end{equation}
\end{corr}
\begin{proof}
As we have seeing before ribbon equivalences between $\cZ(G)$ and $\cZ(Q)$ correspond to algebras $A(U,U,\gamma)$ in
$\cZ(G)\boxtimes\cZ(Q)$ with trivial parents. By applying theorem \ref{pare} we get the conditions of the corollary.
\end{proof}

The next auxiliary result, describing subgroups of direct products, will be used to get a different presentation for
ribbon equivalences.
\begin{lem}\label{subgr}
Subgroups in $G\times Q$ correspond to diagrams of groups
\begin{equation}\label{diag}
\UseTips\xymatrix{G && P && Q\\ & M \ar@{_{(}->}[lu] \ar@{->>}[ru]^i && N \ar@{->>}[lu]_j \ar@{^{(}->}[ru] .}
\end{equation}
The diagram, corresponding to a subgroup $U\subset G\times Q$ has a form:
\begin{equation}\label{diags}
\UseTips\xymatrix{G && P && Q\\ & pr_1(U) \ar@{_{(}->}[lu] \ar@{->>}[ru]^i && pr_2(U) \ar@{->>}[lu]_j \ar@{^{(}->}[ru]
.}
\end{equation}
Conversely, the subgroup, corresponding to a diagram (\ref{diag}), is
$$U = M\times_P N = \{(g,q)\in M\times N|\quad  i(g) = j(q)\}\subset G\times Q.$$
\end{lem}
\begin{proof}
The group $P$ and the surjections in the diagram (\ref{diags}) are defined as follows. First note that, as a subgroup
of $G$, $U\cap(G\times\{e\})$ is a normal subgroup of $pr_1(U)$. This, indeed, follows from the fact that for
$(g,q),(f,e)\in U$
$$(g,q)(f,e)(g,q)^{-1} = (gfg^{-1},e)$$ lies in $U$. Similarly, $U\cap(\{e\}\times Q)$ is a normal subgroup of
$pr_2(U)$. Moreover, there is an isomorphism of quotient groups
\begin{equation}\label{isomm}
pr_1(U)/U\cap(G\times\{e\}) \to pr_2(U)/U\cap(\{e\}\times Q),
\end{equation}
given by the assignment on cosets $g(U\cap(G\times\{e\})) \mapsto q(U\cap(\{e\}\times Q))$ each time $(g,q)$ belongs to
$U$. Thus, in the diagram (\ref{diags}), we can set $P = pr_1(U)/U\cap(G\times\{e\})$ with $i$ being the quotient map
and $j$ being the composition of the quotient map with the inverse of (\ref{isomm}).

The fact that the constructions, described in the lemma, are mutually inverse can be verified directly.
\end{proof}
\begin{rem}
\end{rem}
With the help of lemma \ref{subgr}, the statement of corollary \ref{eqbc} can be reformulated as follows. Equivalences
between $\cZ(G)$ and $\cZ_\beta(Q)$, as ribbon categories, correspond to diagrams
$$\UseTips
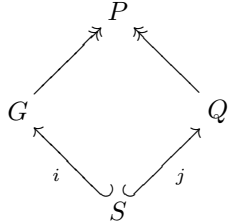
\xymatrix{ & P \\ G \ar@{->>}[ur] && Q  \ar@{->>}[ul] \\ & S \ar@{_{(}->}[ul]^i  \ar@{^{(}->}[ur]_j}$$
with abelian $S$, where the inclusions are normal and such that the actions of $P$ on $S$, induced by the extensions,
coincide; together with the coboundary $\gamma\in C^2(U,k^*)$ on $U = G\times_PQ$ $$d(\gamma) =
(\alpha\times\beta^{-1})|_U,$$ such that $\gamma(\ |\ )$ induces a non-degenerate pairing $i(S)\times j(S)\to k^*$.
\newline
Here we defined $S$ to be $U\cap(G\times\{e\}$ with the obvious inclusion $i$ and with $j$ defined by a choice of
isomorphism $S\to \hat S$ followed by the map $\hat S\to U\cap(\{e\}\times Q)$, induced by the pairing (\ref{pair}).

\begin{rem}
\end{rem}
It is well-known that the category $\cZ(G)$ is equivalent to the monoidal centre $\cZ(\Rep(G))$ of the category
$\Rep(G)$ of (finite-dimensional) representations of $G$. In particular, a monoidal equivalence $\Rep(G)\to \Rep(Q)$
gives rise to an equivalence of ribbon categories $\cZ(G)\to \cZ(Q)$. Monoidal equivalences between categories of
representations of finite groups were described in \cite{da1,da2} (see also \cite{eg}). According to \cite{da2},
monoidal equivalences $\Rep(G)\to \Rep(Q)$ correspond to the following data:  a diagram of groups
$$\UseTips
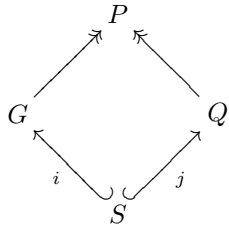
\xymatrix{ & P \\ G \ar@{->>}[ur] && Q  \ar@{->>}[ul] \\ & S \ar@{_{(}->}[ul]^i  \ar@{^{(}->}[ur]_j}$$
with abelian $S$ (such that the actions of $P$ on $S$, induced by the extensions, coincide); together with a
$P$-invariant cohomology class $\gamma\in H^2(S,k^*)^P$ and a homomorphism
$$G\times_P Q\to N(P,S,\gamma) = \{(p,\pi)\in P\times C^1(S,k^*)|\ p(\overline{\gamma})\overline{\gamma}^{-1} = d(\pi)\},$$
(here $\overline{\gamma}\in Z^2(S,k^*)$ is a representative of the class $\gamma$) fitting into the commutative diagram
with exact rows and columns:
\[\begin{array}{ccccc}
 & & P & = & P \\
 & & \uparrow & & \uparrow \\
 S & \to & G\times_PQ & \to & N(P,S,\gamma) \\
 \| & & \uparrow & & \uparrow \\
 S & \to & S\times S & \to & \hat S ,
\end{array}\]
where $S\to S\times S$ is the diagonal embedding and $S\times S\to\hat S$ is a skew-diagonal projection given by the
pairing $\overline{\gamma}(\ |\ ):S\to\hat S$.
\newline
The group-theoretical data of the corresponding ribbon equivalence $\cZ(G)\to \cZ(Q)$ is given by the 2-class
$\tilde\gamma\in H^2(G\times_PQ,k^*)$, which can be constructed using the short exact sequence $S\to G\times_PQ\to
N(P,S,\gamma)$. The details will appear elsewhere.
\begin{rem}
\end{rem}
Since the category $\cZ(G)$ is isomorphic to the monoidal centre $\cZ(\cC(G,\alpha))$, any monoidal equivalence
$\cC(G,\alpha)\to\cC(Q,\beta)$ gives a ribbon equivalence $\cZ(G)\to\cZ_\beta(Q)$. Monoidal equivalences
$\cC(G,\alpha)\to\cC(Q,\beta)$ correspond to isomorphisms $\phi:G\to Q$ together with a coboundary $\gamma\in
C^2(G,k^*)$ $d(\gamma) = \alpha\phi^*(\beta)$. It is straightforward to see, that the corresponding subgroup $U\subset
G\times_PQ$ is the graph of $\phi$ $U = \{(g,\phi(g))|\ g\in G\}$ and that the coboundary $\tilde\gamma\in C^2(U,k^*)$
is given by $\tilde\gamma((f,\phi(f)),(g,\phi(g))) = \gamma(f,g)$.

\section{Modular invariants for group-theoretical modular data}\label{migrth}

\subsection{Characters of commutative algebras and their local modules}

\begin{prop}\label{transcha}
The map $K_0(\cZ(H))\cong K_0(\cZ(G)_{k(H)})\to K_0(\cZ(G))$, induced by the transfer $\cZ(H)\cong
\cZ(G)_{k(H)}\to\cZ(G)$, sends a character $\chi\in K_0(\cZ(H))$ into
$$\overline\chi(f,g) = \frac{1}{|H|}\sum_{x\in G,\\ xfx^{-1},xgx^{-1}\in H}\chi(xfx^{-1},xgx^{-1}).$$
\end{prop}
\begin{proof}
The proof is completely analogous to the prove of the induction formula in character theory (see for example
\cite{is}). By the definition, the character $\overline\chi(f,g)$ is the trace $tr_{(k(G)\otimes_H U)_f}(g)$ of $g$
acting on the graded component $(k(G)\otimes_H U)_f$. By the definition of the transfer, the graded component
$(k(G)\otimes_H U)_f$ coincides with $\oplus_x p_x\otimes U_{x^{-1}fx}$, where the sum is taken over cosets $\{x:
x^{-1}fx\in H\}/H$ (with respect to the $H$-action on the set $\{x: x^{-1}fx\in H\}$ by left multiplications). So that
$$tr_{(k(G)\otimes_H U)_f}(g) = \sum_x tr_{p_x\otimes U_{x^{-1}fx}}(g).$$ Note that $g$ preserves $p_x\otimes
U_{x^{-1}fx}$ iff $x^{-1}gx$ is in $H$:
$$g(p_x\otimes U_{x^{-1}fx}) = p_{gx}\otimes U_{x^{-1}fx} = p_{xx^{-1}gx}\otimes U_{x^{-1}fx} = p_x\otimes x^{-1}gx(U_{x^{-1}fx}),$$
and that, in this case, the restriction of $g$ to $p_x\otimes U_{x^{-1}fx}$ coincides with $p_x\otimes x^{-1}gx$. Thus
$$tr_{p_x\otimes U_{x^{-1}fx}}(g) = \chi(xfx^{-1},xgx^{-1}).$$
\end{proof}

\begin{corr}
The character of a trivialising algebra $A(H,\gamma)$ has the form:
\begin{equation}\label{chalg}
\chi_{A(H,\gamma)}(f,g) = \frac{1}{|H|}\sum_{x\in G,\\ xfx^{-1},xgx^{-1}\in H}\gamma(xfx^{-1}|xgx^{-1}).
\end{equation}
\end{corr}
\begin{proof}
By the definition the algebra $A(H,\gamma)$ is the image of the algebra $k[H,\gamma]$ under the transfer $\cZ(H)\cong
\cZ(G)_{k(H)}\to\cZ(G)$. Thus the corollary follows from proposition \ref{transcha} and the fact that the character of
$k[F,\gamma]$ is $\chi(x,y)=\gamma(x|y)$, which can be checked directly. Indeed, $x$-graded component of $k[F,\gamma]$
is spanned (over $k$) by $e_x$, with the action of $y$ on it
$$y(e_x) = \varepsilon_x(y)e_{yxy^{-1}} = \frac{\gamma(x,y)}{\gamma(xyx^{-1},x)}e_x.$$ Thus, for commuting $x,y$, we
have that
$$tr_{k[F,\gamma]_x}(y) = \frac{\gamma(x,y)}{\gamma(xyx^{-1},x)} = \gamma(x|y).$$
\end{proof}

We finish this section with examples of trivialising algebras in $\cZ(G)$ with the same character (the same class in
$K_0(\cZ(G))$). By the formula \ref{chalg}, to construct such example it is enough to have a finite group $H$ with a
non-trivial 2-class $\gamma\in H^2(H,k^*)$, such that $\gamma_h$ is a trivial character of $C_H(h)$ for any $h\in H$.
We will use the well known correspondence between 2-cohomology and central extensions (see \cite{br}, for example) to
give a group theoretic conditions, which guarantee the existence of such class.
\begin{lem}
Let $\tilde H$ be a finite group with a central cyclic subgroup $Z\subset Z(\tilde H)$, which does not contain
commutators, and such that the extension $Z\cap [\tilde H,\tilde H]\to [\tilde H,\tilde H]\to [\tilde H,\tilde H]/Z\cap
[\tilde H,\tilde H]$ is non-trivial. Let $k$ be an algebraically closed field of characteristic zero. Then the class
$\gamma\in H^2(H,k^*)$, extended from the extension class $\overline\gamma\in H^2(H,Z)$ with respect to an embedding
$Z\to k^*$, is non-trivial and such that $\gamma_h$ is a trivial character of $C_H(h)$ for any $h\in H$. Here $H=\tilde
H/Z$.
\end{lem}
\begin{proof}
We begin by showing that the absence of commutators in $Z$ implies that $\gamma_h$ is a trivial character of $C_H(h)$
for any $h\in H$. Indeed, it is strightforward to see that for any commuting $x,h\in H$ the commutator $[\tilde
h,\tilde x]$ of (any of) their preimages in $\tilde H$ lies in $Z$ and coincides with $\overline\gamma_h(x)$. So if $Z$
does not contain commutators, then $\overline\gamma_h$ is trivial for any $h\in H$, which implies the triviality of
$\gamma_h$.
\newline
Next we show non-triviality of the extended class $\gamma\in H^2(H,k^*)$. We identify $Z$ with the $n$-torsion subgroup
of $k^*$ ($n=|Z|$). It follows from the long exact sequence, corresponding to the coefficient extension $Z\to
k^*\stackrel{\times n}{\to}k^*$, that the kernel of the coefficient extension $H^2(H,Z)\to H^2(H,k^*)$ coincides with
the image of the connecting map $\partial:H^1(H,k^*)\to H^2(H,Z)$, which fits into a diagram
$$\xymatrix{ & & H^2([\tilde H,\tilde H],Z)\\ \ar[r] & H^1(H,k^*) \ar[r]^\partial  & H^2(H,Z) \ar[u] \ar[r] & H^2(H,k^*) \ar[r]& \\
\ar[r] & H^1(H^{ab},k^*) \ar[r] \ar[u] & H^2(H^{ab},Z) \ar[u] \ar[r] & }$$ Here $H^{ab} = H/[H,H]$. The commutative
square implies that the image of $\partial$ coincides with the image of $H^2(H^{ab},Z)\to H^2(H,Z)$. So if the image of
$\gamma\in H^2(H,Z)$ in $H^2([\tilde H,\tilde H],Z)$ is non-trivial, then $\gamma$ can not be in the image of
$H^2(H^{ab},Z)\to H^2(H,Z)$ and thus can not be killed by $H^2(H,Z)\to H^2(H,k^*)$.
\end{proof}
\begin{exam}
\end{exam}
Let $p\geq 5$ be a prime and $\tilde H$ be the free meta-abelian group of period $p^2$, generated by $x_1,x_2,x_3,x_4$.
Let $Z$ be the central subgroup, generated by $([x_1,x_2][x_3,x_4])^p$. Note that $V=\tilde H^{ab}$ is the free abelian
group of period $p^2$, with four generators $e_1,e_2,e_3,e_4$, and $[\tilde H,\tilde H]$ can be identified with the
exterior square $\Lambda^2V$. In this presentation the commutator pairing correspond to the wedge product $V\times V\to
\Lambda^2V$ so the set of commutators in $[\tilde H,\tilde H]$ correspond to the Pl\"ucker quadric $\{x\in \Lambda^2V,\
x\wedge x=0\}\subset \Lambda^2V$. The element $v=p(e_1\wedge e_2+e_3\wedge e_4)$ is not on the quadric, which shows
that $Z$ does not contain (non-trivial) commutators. The inclusion $\langle v\rangle\to \Lambda^2V$ does not split,
which implies that the extension $Z\cap [\tilde H,\tilde H]\to [\tilde H,\tilde H]\to [\tilde H,\tilde H]/Z\cap [\tilde
H,\tilde H]$ is non-trivial.
\newline
Applying the lemma we get a desired example.

\subsection{$S_3$ modular data and modular invariants}

Recall that $H^2(H,k^*)$ is trivial for any subgroup $H$ of $S_3$ (including $S_3$ itself). The classes of simple
objects are labeled by
$$(e,\xi_0),(e,\xi_1),(e,\xi_2),((123),\pi_0),((123),\pi_1),((123),\pi_0),((12),\psi_0),((12),\psi_1).$$
Here $\xi_i\in Irr(S_3),\ \pi_i\in Irr(C_3),\ \psi_i\in Irr(C_2).$

The $S$- and $T$-matrices have the following form:
$$S = \frac{1}{6}\left(\begin{array}{rrrrrrrr}
1 & 1 & 2 & 2 & 2 & 2 & 3 & 3\\ 1 & 1 & 2 & 2 & 2 & 2 & -3 & -3\\ 2 & 2 & 4 & -2 & -2 & -2 & 0 & 0\\ 2 & 2 & -2 & 4 &
-2 & -2 & 0 & 0\\ 2 & 2 & -2 & -2 & -2 & 4 & 0 & 0\\ 2 & 2 & -2 & -2 & 4 & -2 & 0 & 0\\ 3 & -3 & 0 & 0 & 0 & 0 & 3 &
-3\\ 3 & -3 & 0 & 0 & 0 & 0 & -3 & 3\\
\end{array}\right)$$
$$T = diag(1,1,1,1,\omega,\omega^{-1},1,-1),\quad 1+\omega+\omega^2=0.$$

In the table below we list all indecomposable commutative separable algebras in $\cZ(S_3)$ together with the characters
of their simple local modules (the first character is the character of the algebra itself):

$$\begin{array}{c|c|l} H\vartriangleright F & H/F & {\cZ(S_3)}^{loc}_{A(H,F)} \\ \hline
S_3\vartriangleright S_3 & \{e\} & \chi_0+\chi_3+\chi_6 \\ A_3\vartriangleright A_3 & \{e\} & \chi_0+\chi_1+2\chi_3 \\
C_2\vartriangleright C_2 & \{e\} & \chi_0+\chi_2+\chi_6 \\ \{e\}\vartriangleright \{e\} & \{e\} & \chi_0+\chi_1+2\chi_2
\\ S_3\vartriangleright A_3 & C_2 & \chi_0+\chi_3, \chi_1+\chi_3, \chi_6, \chi_7 \\ C_2\vartriangleright \{e\}  & C_2 &
\chi_0+\chi_2, \chi_1+\chi_2, \chi_6, \chi_7 \\ A_3\vartriangleright \{e\} & A_3 & \chi_0+\chi_1, \chi_2, \chi_3,
\chi_4, \chi_5 \\ S_3\vartriangleright \{e\} & S_3 & \chi_0, \chi_1, \chi_2, \chi_3, \chi_4, \chi_5, \chi_6, \chi_7
\end{array}$$
Note, that in the case $A_3\vartriangleright \{e\}$, for each $i=2,...,5$ there are two different simple local modules
with the character $\chi_i$.

According to lemma \ref{subgr}, there are 22 conjugacy classes of subgroups in $S_3\times S_3$. Sixteen of them have a
form $A\times B$ for $A,B\subset S_3$ and correspond to diagrams (we omit the embeddings into $S_3$)
$$\UseTips
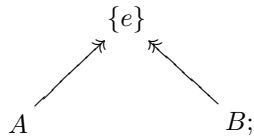
\xymatrix{ & \{e\} \\ A \ar@{->>}[ru] & & B ;\ar@{->>}[ul]  }$$ four have the form
$$\delta(C_2),\ \delta(C_2)(\{e\}\times A_3),\ \delta(C_2)(A_3\times \{e\}) ,\ \delta(C_2)(A_3 \times A_3)$$
and correspond to diagrams
$$\UseTips
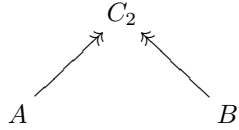
\xymatrix{ & C_2 \\ A \ar@{->>}[ru] & & B \ar@{->>}[ul]}$$
with $A,B = C_2$ or $S_3$; and two remaining are $\delta(A_3), \delta(S_3)$, corresponding to the diagrams:
$$\UseTips
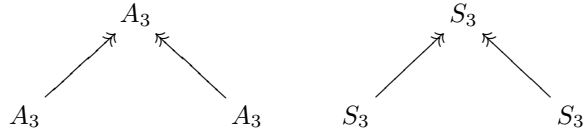
\xymatrix{ & A_3 &&& S_3\\ A_3 \ar@{->>}[ru] & & A_3 \ar@{->>}[ul] & S_3 \ar@{->>}[ru] & & S_3 \ar@{->>}[ul]}$$
respectively. Only six of them have non-trivial cohomology $H^2(U,k^*)$:
$$C_2\times C_2, S_3\times S_3, C_2\times S_3, S_3\times C_2, A_3\times A_3, \delta(C_2)(A_3\times A_3).$$
For the first four $H^2(U,k^*)$ is cyclic of order 2 and for two remaining it is cyclic of order 3. In the last two
cases the conjugation action of the normaliser in $S_3\times S_3$ permutes two non-trivial cohomology classes. Thus, in
all cases, there is just one (up to conjugation) non-trivial cohomology class, which, somewhat loosely, will be denoted
$\gamma$.

Maximal commutative algebras in $\cZ(S_3)\boxtimes \cZ(S_3) \simeq \cZ(S_3\times S_3)$ are depicted as edges of the
following four graphs:

$$\UseTips
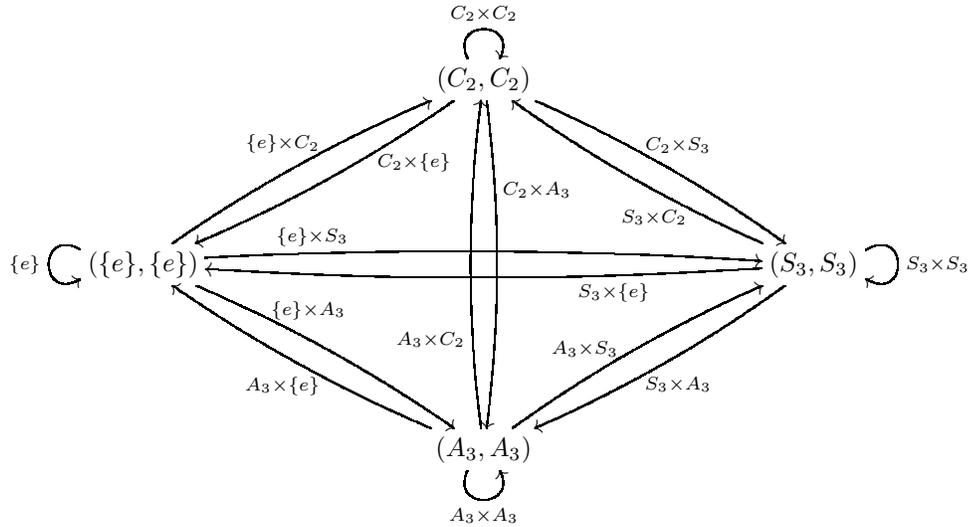
\xymatrix{ &&& (C_2,C_2) \ar@(ul,ur)[]^{C_2\times C_2} \ar@/^5pt/[llldd]^(0.3){C_2\times \{e\}} \ar@/^5pt/[rrrdd]^{C_2\times S_3}
\ar@/^5pt/[dddd]^(0.3){C_2\times A_3} \\ \\ (\{e\},\{e\}) \ar@/^5pt/[rrruu]^{\{e\}\times C_2}
\ar@/^5pt/[rrrdd]^(0.4){\{e\}\times A_3}
 \POS!L(.7)="a",\ar@(ul,dl)"a"_{\{e\}}
\ar@/^5pt/[rrrrrr]^(0.3){\{e\}\times S_3}  &&&&&& (S_3,S_3) \ar@/^5pt/[llllll]^(0.3){S_3\times\{e\}}
\ar@/^5pt/[ddlll]^{S_3\times A_3} \ar@/^5pt/[uulll]^(0.4){S_3\times C_2} \POS!R(.65)="b",\ar@(ur,dr)"b"^{S_3\times S_3}
\\ \\
&&& (A_3,A_3) \ar@/^5pt/[uuuu]^(0.3){A_3\times C_2} \ar@/^5pt/[uurrr]^(0.4){A_3\times S_3} \ar@/^5pt/[uulll]^{A_3\times
\{e\}} \ar@(dl,dr)[]_{A_3\times A_3} }$$

$$\UseTips
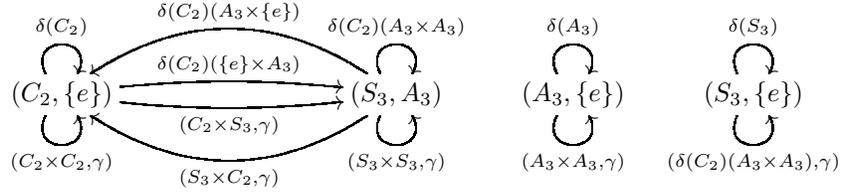
\xymatrix{(C_2,\{e\}) \ar@(ul,ur)[]^{\delta(C_2)} \ar@(dl,dr)_{(C_2\times C_2,\gamma)}
\ar@/^5pt/[rrr]^{\delta(C_2)(\{e\}\times A_3)} \ar@/_5pt/[rrr]_{(C_2\times S_3,\gamma)} &&& (S_3,A_3)
\ar@/_25pt/[lll]_{\delta(C_2)(A_3\times \{e\})} \ar@/^25pt/[lll]^{(S_3\times C_2,\gamma)}
\ar@(ul,ur)[]^{\delta(C_2)(A_3\times A_3)} \ar@(dl,dr)_{(S_3\times S_3,\gamma)} & (A_3,\{e\})
\ar@(ul,ur)[]^{\delta(A_3)} \ar@(dl,dr)_{(A_3\times A_3,\gamma)} & (S_3,\{e\}) \ar@(ul,ur)[]^{\delta(S_3)}
\ar@(dl,dr)_{(\delta(C_2)(A_3\times A_3),\gamma)}}$$ Vertices are labeled by the (conjugacy classes of) pairs of
subgroups $F\vartriangleleft H\subset S_3$, which correspond to indecomposable commutative separable algebras in
$\cZ(S_3)$; edges are labeled by (conjugacy classes of) $(H,\gamma)$, where $H\subset S_3\times S_3$ and $\gamma\in
H^2(H,k^*)$ is a cohomology class (omitted if trivial), which correspond to maximal indecomposable commutative
separable algebras in $\cZ(S_3\times S_3)$ . An edge goes from $A$ to $B$ if for the corresponding algebra $C\in
\cZ(S_3\times S_3)$ $$Hom_{1\boxtimes \cZ(S_3)}(1,C) = A,\quad Hom_{\cZ(S_3)\boxtimes 1}(1,C) = B.$$

The table below contains characters of maximal indecomposable commutative separable algebras in $\cZ(S_3\times S_3)$,
written in the (traditional) form of a partition function, ordered by the rank of the corresponding modular invariant.

$$\begin{array}{l|l} A(H,\gamma) & Z \\ \hline
A(S_3\times S_3) &  |\chi_0+\chi_3+\chi_6|^2 \\ A(S_3\times A_3) &  (\chi_0+\chi_3+\chi_6)(\chi_0+\chi_1+2\chi_3)^* \\
A(S_3\times C_2) &  (\chi_0+\chi_3+\chi_6)(\chi_0+\chi_2+\chi_6)^* \\ A(S_3\times\{e\}) &
(\chi_0+\chi_3+\chi_6)(\chi_0+\chi_1+2\chi_2)^* \\ A(A_3\times A_3) &  |\chi_0+\chi_1+2\chi_3|^2 \\ A(A_3\times S_3) &
(\chi_0+\chi_1+2\chi_3)(\chi_0+\chi_3+\chi_6)^* \\ A(A_3\times C_2) &  (\chi_0+\chi_1+2\chi_3)(\chi_0+\chi_2+\chi_6)^*
\\ A(A_3\times\{e\}) &  (\chi_0+\chi_1+2\chi_3)(\chi_0+\chi_1+2\chi_2)^* \\ A(C_2\times C_2) &
|\chi_0+\chi_2+\chi_6|^2 \\ A(C_2\times S_3) &  (\chi_0+\chi_2+\chi_6)(\chi_0+\chi_3+\chi_6)^* \\ A(C_2\times A_3) &
(\chi_0+\chi_2+\chi_6)(\chi_0+\chi_1+2\chi_3)^* \\ A(C_2\times\{e\}) &  (\chi_0+\chi_2+\chi_6)(\chi_0+\chi_1+2\chi_2)^*
\\ A(\{e\}\times \{e\}) &  |\chi_0+\chi_1+2\chi_2|^2 \\ A(\{e\}\times S_3) &
(\chi_0+\chi_1+2\chi_2)(\chi_0+\chi_3+\chi_6)^* \\ A(\{e\}\times A_3) &
(\chi_0+\chi_1+2\chi_2)(\chi_0+\chi_1+2\chi_3)^* \\ A(\{e\}\times C_2) &
(\chi_0+\chi_1+2\chi_2)(\chi_0+\chi_2+\chi_6)^* \\ A(\delta(C_2)(A_3\times A_3)) &
|\chi_0+\chi_3|^2+|\chi_1+\chi_3|^2+|\chi_6|^2+|\chi_7|^2 \\ A(S_3\times S_3,\gamma) &
|\chi_0+\chi_3|^2+(\chi_1+\chi_3)\chi_6^*+\chi_6(\chi_1+\chi_3)^*+|\chi_7|^2 \\ A(\delta(C_2)(A_3\times\{e\})) &
(\chi_0+\chi_3)(\chi_0+\chi_2)^*+(\chi_1+\chi_3)(\chi_1+\chi_2)^*+|\chi_6|^2+|\chi_7|^2 \\ A(S_3\times C_2,\gamma) &
(\chi_0+\chi_3)(\chi_0+\chi_2)^*+(\chi_1+\chi_3)\chi_6^*+\chi_6(\chi_1+\chi_2)^*+|\chi_7|^2 \\ A(\delta(C_2)) &
|\chi_0+\chi_2|^2+|\chi_1+\chi_2|^2+|\chi_6|^2+|\chi_7|^2 \\ A(S_3\times S_3,\gamma) &
|\chi_0+\chi_2|^2+(\chi_1+\chi_2)\chi_6^*+\chi_6(\chi_1+\chi_2)^*+|\chi_7|^2 \\ A(\delta(C_2)(\{e\}\times A_3)) &
(\chi_0+\chi_2)(\chi_0+\chi_3)^*+(\chi_1+\chi_2)(\chi_1+\chi_3)^*+|\chi_6|^2+|\chi_7|^2 \\ A(C_2\times S_3,\gamma) &
(\chi_0+\chi_2)(\chi_0+\chi_3)^*+(\chi_1+\chi_2)\chi_6^*+\chi_6(\chi_1+\chi_3)^*+|\chi_7|^2 \\ A(\delta(A_3)) &
|\chi_0+\chi_1|^2+2|\chi_2|^2+2|\chi_3|^2+2|\chi_4|^2+2|\chi_5|^2 \\ A(A_3\times A_3,\gamma) &
|\chi_0+\chi_1|^2+2\chi_2\chi_3^*+2\chi_3\chi_2^*+|\chi_4|^2+|\chi_5|^2 \\ A(\delta(S_3)) &
|\chi_0|^2+|\chi_1|^2+|\chi_2|^2+|\chi_3|^2+|\chi_4|^2+|\chi_5|^2+|\chi_6|^2+|\chi_7|^2\\ A(\delta(S_3)(A_3\times
A_3),\gamma) &  |\chi_0|^2+|\chi_1|^2+\chi_2\chi_3^*+\chi_3\chi_2^*+|\chi_4|^2+|\chi_5|^2+|\chi_6|^2+|\chi_7|^2
\end{array}$$

\section{Concluding remarks}
It is known that modular categories give rise to 3-dimensional topological field theories (see \cite{tu,bk}). In
particular, group-theoretical modular categories correspond to Dijkgraaf-Witten topological field theories
\cite{dw,fr}. Equivalent modular categories give rise to equivalent topological field theories and, in particular, to
the same invariants of closed 3-manifolds. It follows from the results of \cite{dw,fr} that the invariant of a
3-manifold $M$, defined by the modular category $\cZ(G)$, has the form
\begin{equation}\label{i3m}
Z_{G}(M) = \frac{|Hom(\pi_1(M),G)|}{|G|}.
\end{equation}
It follows from corollary \ref{eqbc} that for $G$ and $Q$, satisfying the conditions of the corollary, the invariants
coincide $$Z_{G}(M) = Z_{Q}(M)$$ for all closed 3-manifolds $M$. In particular, the number of homomorphisms
$\pi_1(M)\to G$ is equal to the number of homomorphisms $\pi_1(M)\to Q$, for $G$ and $Q$ satisfying the conditions of
corollary \ref{eqbc}. In other words, fundamental groups of 3-manifolds {\em do not feel the difference} between such
$G$ and $Q$.

\end{document}